\documentclass[final, 3p, times]{elsarticle}
\usepackage{amssymb}
\usepackage{multirow}
\usepackage{setspace}
\usepackage{amsthm}
\usepackage{amsmath}
\usepackage[centerlast]{subfigure}
\usepackage{float}
\usepackage{stfloats} 
\usepackage{algorithm}
\usepackage{algorithmic}
\usepackage{wrapfig}

\usepackage{amssymb}

\journal{Structural Safety}

\begin{document}

\begin{frontmatter}



\title{Reliability sensitivity analysis based on probability distribution perturbation with application to $\mbox{CO}_{2}$ storage}

\author[UT,IFP]{Ekaterina Sergienko}
\author[EDF,INRIA]{Paul Lema\^ itre}
\author[EDF]{Aur\'elie Arnaud}
\author[IFP]{Daniel Busby}
\author[UT]{Fabrice Gamboa}

\address[UT]{Universit\'e Paul Sabatier, IMT-EPS, 118 Route de Narbonne, 31062, Toulouse, France}
\address[IFP]{IFP Energies Nouvelles, 1-4 Avenue de Bois-Pr\'eau, 92582, Rueil-Malmaison, France}
\address[EDF]{EDF R\&D, 6 Quai Watier, 78401, Chatou, France}
\address[INRIA]{INRIA Sud-Ouest, 351 Cours de la Lib\'eration, 33405, Talence, France}

\begin{abstract}
The objective of reliability sensitivity analysis is  to determine input variables that mostly contribute  to the variability of the
failure probability. In this paper, we study a recently introduced method for the reliability sensitivity analysis based on a perturbation of the original probability distribution of the input variables. The objective is to determine the most influential input variables and to analyze their impact on the failure probability. We propose a moment independent sensitivity measure that is based on a perturbation of the original probability density independently for each input variable. The variables providing the highest variation of the original failure probability are settled to be more influential. These variables will need a proper characterization in terms of uncertainty. The method is intended to work in applications involving a computationally expensive simulation code for evaluating the failure probability such as the $\mbox{CO}_{2}$ storage risk analysis. An application of the method to a synthetic $\mbox{CO}_{2}$ storage case study is provided together with some analytical examples. 
\end{abstract}

\end{frontmatter}
\section{Introduction}
Carbon Capture and Storage (CCS) stands for the collection of $\mbox{CO}_{2}$ from industrial sources and its injection into deep geological formations for a  permanent storage. There are  three possible sites for injection: unmined coalbed formations, saline aquifers  and depleted oil and gas reservoirs  \citep{polson2009}. Nevertheless, the following principal environmental question arises: what is the probability that $\mbox{CO}_{2}$ will remain underground for hundreds to thousands of years after its capture and injection into a storage formation?  \par
The primary risk of  $\mbox{CO}_{2}$ geological storage is unintended gas leakage from the storage reservoir \citep{leguen2008, bowden}. In this work, we focus on a leakage from the storage formation through a fault or a fracture. This can happen when the reservoir pressure is higher than the caprock fracture pressure.  Numerical modelling and simulation has become an integral component of $\mbox{CO}_{2}$ storage assessment and monitoring. The reservoir simulation models are constructed based on the reservoir and production data. They are  used to predict and analyze $\mbox{CO}_{2}$ plume distribution and the reservoir pressure development during the injection and the storage periods. With the help of the numerical simulation models it is possible to forecast $\mbox{CO}_{2}$ storage performance and to evaluate the risks of a possible leakage.  
\par
Risk and uncertainty analysis has been recognized as a principal part of safety and risk assessment.  Generally speaking, when the main sources of uncertainty have been identified, Uncertainty Analysis (UA) is focused on quantifying the uncertainty in the model output resulting from uncertainty in the model inputs. At the same time, Sensitivity Analysis (SA) aims to identify  the contributions of each model uncertain input to the variability of the model output. Uncertainty analysis may be equally performed  to assess the reliability of the system. A typical example of a failure probability estimation in the $\mbox{CO}_{2}$ storage risk analysis is the estimation of the probability of exceeding the caprock fracturing pressure during the $\mbox{CO}_{2}$ injection phase. If we denote $\mbox{P}_{reservoir}$ the reservoir pressure and $\mbox{P}_{fracture}$ the caprock fracturing pressure, then we consider the following failure probability:
\[ p_{f}=\mbox{P}(\mbox{P}_{reservoir}\geq \mbox{P}_{fracture}).\]
In practice, uncertainty and sensitivity analysis  require a large number of reservoir simulator runs  to explore all the input variables space.
However, higher accuracy of a simulator usually results in a higher simulation time.  One simulator run can take from few minutes up to several hours or even days. Therefore, when the simulation time becomes too high, uncertainty analysis may become unfeasible. For this reason, in this work, in order to estimate the failure probability $p_{f}$ for the expensive reservoir simulator, we use Gaussian Process (GP) response surface model, also known as \textsl{ kriging} \cite{Sacks1989, Welch1992, Santner2003}. The kriging method was originally introduced in the field of geostatistics by Krige in the 1950's \citep{krig} and formalized in 1960's by Matheron \cite{Matheron1963}. In \cite{Sacks1989} Sacks et al. proposed the statistical approach to uncertainty analysis of complex computer codes referred to as the Design  of Computer Experiments. In a nutshell, the approach consists in building an approximation of the reservoir simulator input/output relationship starting from a set of simulation runs at a carefully chosen input variables configurations referred to as experimental design or training set. The obtained response surface model can then be used to predict the model output for a new non simulated input with a negligible computational time.  Therefore, uncertainty and sensitivity analysis become affordable.  
\par
One of the most challenging problems in risk analysis is to identify the failure region and to compute the failure probability. However, the failure probability will usually depend strongly on the probability distribution of the input variables. Reliability Sensitivity Analysis can help in understanding the relationship between each input variable uncertainty and the failure probability. The problem is to identify the set of input variables that need to be well characterized in terms of uncertainty distribution.  
\par For the time being, numerous approximation and simulation methods are available for estimating the failure probability such as First / Second Order Reliability Methods, Monte Carlo sampling, importance sampling, directional sampling, subset simulation, etc. \citep{lemaire2009, melchers89,ditlevsen88, beck2001}. However, there are few sensitivity analysis methods developed for the failure probability analysis. The widely used methods for sensitivity analysis are based on a variance decomposition of the output.  Given the probability distribution of the input variables, Sobol indices are expressed by the ratio of the variance due to a given input on the total output variance \cite{sobol1990}.  Knowing the probability distribution of the input variables and the output, we can so define sensitivity indices for each of the input variables. Variance based indices are usually interesting for measuring the input/output sensitivity, however they can be poorly relevant to our problem of evaluating the impact on failure probability.   
\par
There have been few attempts to develop sensitivity analysis methods well suited for reliability analysis. First, as complementary results of the First Order Reliability Method, sensitivity to the distribution of the input variables can be obtained. Sensitivity is expressed as the partial derivative of the reliability index $\boldsymbol{\beta}$ \citep{lemaire2009}. Another approach was proposed by Morio \cite{morio2011influence}. Therein, the author uses the variance decomposition and Sobol' sensitivity indices to study the rate of change in the failure probability due to the changes of the input distribution density parameters. Borgonovo et al. \citep{borgonovo2011moment} suggested some moment independent importance measure in the reliability analysis. This measure does not involve the variance. For a fixed variable $x_{i}$ it quantifies the effect of knowing $x_{i}$ by computing the $L_{1}$ norm between the unconditional joint density $f_{\bar{x}}(\cdot)$ and the conditional density $ f_{\bar{x}|x_{i}}(\cdot)$.
\par
In this paper, we study a moment independent approach for sensitivity analysis of a failure probability \citep{paul2012}. The influence of the input variables on the failure probability is obtained by perturbing the prior probability density function $f_{\bar{x}}(\bar{x})$. In particular, we estimate the effect of the perturbation on the value of the failure probability $p_{f}$. Here, we propose to distinguish distributions classes by their supports. For the case of a bounded support, such as the uniform or the triangular distributions, the main source of uncertainty is about the boundaries of the support. On the other hand, in the case of infinite support, such as normal or log-normal distributions, the main source of uncertainty comes from the distribution parameters, such as mean and variance. The estimation method has the advantage of being very efficient in terms of number of simulator calls. In order to estimate the sensitivity indices for all the input variables, the performance function is evaluated only once on a Monte Carlo sample used to estimate the reference failure probability $p_{f}$.
\par
Our paper is organized as follows. First, we introduce some density perturbations for different families of distributions. Later, we introduce the technique to compute a perturbed failure probability using the same Monte Carlo sample. This is based on an inverse importance sampling technique \citep{hesterberg1996estimates}. Finally, we present the formulation for the moment independent sensitivity indices and demonstrate its applicability on an analytical and a $\mbox{CO}_{2}$ storage reservoir case examples. 

\section{Density perturbation influence to failure probability}
Let us denote by $g(\bar{x})$ the performance function of the system, $\bar{x}\in\Omega\subset\mathbb{R}^{d}$ is a set of independent input variables with the joint density $f_{\bar{x}}(\bar{x})=\prod_{i=1}^{d}f_{x_{i}}(x_{i})$. The failure probability is expressed as:
\[
p_{f}=\mbox{P}(g(\bar{x})\leq0)=\mathbb{E}_{f_{\bar{x}}}\left[\mathbf{I}_{g(\bar{x})\leq0}\right]=\int_{\Omega_{f}}f_{\bar{x}}(\bar{x})d\bar{x},
\]
where $\Omega_{f}=\{\bar{x}\in\Omega:g(\bar{x})\leq0\}$ is the failure region.
\par
In general, the distribution density $f_{\bar{x}}(\cdot)$ is provided by experts on the basis of some indirect measurements or some limited observation data \citep{lemaitreanalyse}. Here, we study how a perturbation of the original probability density $f_{\bar{x}}(\cdot)$ affects the failure probability of the system $p_{f}$. We assume that the input variables $x_{i}, i=1,\ldots,d$ are independent random variables with marginal densities $f_{x_{i}}$. So that, $f_{\bar{x}}(\bar{x})=\prod_{i=1}^{d}f_{x_{i}}(x_{i})$. 
\par
This work is inspired by previous work \cite{lemaitreanalyse, paul2012}. Originally the method was applicable only to the normal distributions. Here, the method is extended to the more general case of an exponential family. The objective of this study is to estimate the influence of the input random variable from the vector $\bar{x}$ to the failure probability of the system $p_{f}$. The proposed approach consists in perturbing the original density for a given fixed variable $x_{i}$ while keeping constant the probability density functions for all the other variables $\bar{x}_{-i}=(x_{1},\ldots,x_{i-1},x_{i+1},\ldots,x_{d})$. Then, a new value for the failure probability is computed. If this new value $p_{i\delta}$ differs significantly from the reference value $p_{f}$, it means that this selected input variable $x_{i}$ is influential. Conversely, if the new failure probability $p_{i\delta}$ is close to $p_{f}$, then the input $x_{i}$ has low influence on the failure probability.

\subsection{Density perturbation for an exponential family}
To define the probability density perturbation, first recall the definition of an exponential family. An exponential family is a set of distribution having density function that can be expressed in the form of:
\begin{equation}
\label{expfam}
f_{x}(x|\boldsymbol{\theta})=h(x)\exp\left(\eta(\boldsymbol{\theta})\mbox{T}(x)-A(\boldsymbol{\theta})\right),
\end{equation}
where $\boldsymbol{\theta}$ is a vector of distribution parameters, $\eta(\boldsymbol{\theta})$ is a natural parameter, $\mbox{T}(x)$ is a vector of sufficient statistics, $h(x)$ is an underlying weight function and $A(\boldsymbol{\theta})$ is the cumulant generating function \citep{MR489333}. The cumulant generating function ensures that the distribution integrates to one, i.e.:
\[ A(\boldsymbol{\theta})=\log\int h(x)\exp(\eta(\boldsymbol{\theta})\mbox{T}(x))\mu(dx),\]
where $\mu$ is the reference measure (for example Lebesgue measure). It could be continuous or discrete. In this paper, we mainly consider continuous measure. However, all the calculations are valid for a discrete measure as well. Exponential family contains most of the standard discrete and continuous distributions that we use for practical modelling, such as the normal, Poisson, Binomial, exponential, Gamma, multivariate normal, etc.
\par
To define the density perturbation for this family of distributions, we use some ideas coming from information theory \citep{MR1122806}. Kullback-Leibler (KL) divergence is used to measure the magnitude of a perturbation. KL divergence quantifies the "closeness" of two probability distribution $P$ and $Q$. Suppose that $P$ and $Q$ are continuous probability distributions with densities $p(\bar{x})$ and $q(\bar{x})$ (with respect to Lebesgue measure). Then, the KL divergence between $P$ and $Q$ is given by:
\begin{equation}
\mbox{D}_{\mbox{KL}}(P,Q)=\int_{-\infty}^{\infty}p(\bar{x})\log\frac{p(\bar{x})}{q(\bar{x})}d\bar{x}
\end{equation}
For $\delta>0$, let us denote for the variable $x_{i}$ the new (perturbed) density as $f_{x_{i\tau}}(\cdot)$,. We select $f_{x_{i\tau}}(\cdot)$ in such way that:
\begin{equation}
\label{kl}
\mbox{D}_{\mbox{KL}}(f_{x_{i\tau}},f_{x_{i}})=\delta.
\end{equation}
Possible values of the perturbation $\delta$ may be restricted by some inequalities on Kullback-Leibler divergence \citep{dragomir2000}. If we define the function $r(\bar{x}):=\frac{p(\bar{x})}{q(\bar{x})}, \left( x\in\Omega\subset\mathbb{R}^{d}\right)$ and assume that $0<r<r(\bar{x})<R$ for all $\bar{x}\in\Omega\subset\mathbb{R}^{d}$. Then, according to \citep{dragomir2000}, we have:
\begin{eqnarray*}
\mbox{D}_{\mbox{KL}}(p,q)&\geq& 0\\
\mbox{D}_{\mbox{KL}}(p,q)&\leq& \frac{(R-r)^{2}}{4rR}=\delta_{max}
\end{eqnarray*}
According to these inequalities, we choose $\delta\in[0,\delta_{max}]$, where $\delta_{max}= \frac{(R-r)^{2}}{4rR}$ can be computed precisely.
\par
Let us consider the original density $f_{x_{i}}(x)=h(x)\exp(\eta(\boldsymbol{\theta})\mbox{T}(x)-A(\boldsymbol{\theta}))$ from exponential family. In order to stay in the same family, we propose to restrict the choice of possible perturbations among the following class of densities:
\begin{equation}
\label{dkl}
f_{x_{i\tau}}(x)=\exp(\tau \mbox{T}(x)-\psi(\tau))f_{x_{i}}(x)=h(x)\exp\left(\mbox{T}(x)(\eta(\boldsymbol{\theta})-\tau)-(A(\boldsymbol{\theta})+\psi(\tau))\right),
\end{equation}
Here, $\tau$ is a constant depending on $\delta$ (it is chosen under the condition \eqref{kl}). The function $\psi(\tau)$ is a normalization function and it may be expressed as:
\[
\psi(\tau)=\log\left[\int_{-\infty}^{\infty}\exp(\tau \mbox{T}(t))f_{x_{i}}(t)dt\right].\]
It is the cumulant generating function for the perturbed probability  distribution $f_{x_{i\tau}}(x)$. Moreover, if $\mu$ and $\sigma^{2}$ are the mean and the variance of the original probability distribution, then:
\begin{equation}
\label{psiprops}
\begin{array}{c}
\psi(0)=0\\
\psi'(0)=\mu \\
\psi''(0)=\sigma^{2}
\end{array}
\end{equation}\noindent
We aim to perturb $f_{x_{i}}$ in such a way that the KL divergence between the original density $f_{xi}$ and the perturbed density $f_{x_{i\tau}}$ is equal to $\delta$. Notice that:
\begin{eqnarray}
\mbox{D}_{\mbox{KL}}(f_{x_{i\tau}},f_{x_{i}})&=&\int_{-\infty}^{\infty}f_{x_{i\tau}}(t)\log\frac{f_{x_{i\tau}}(t)}{f_{x_{i}}(t)}dt =\int_{-\infty}^{\infty}f_{x_{i\tau}}(t)\left(\tau \mbox{T}(t) - \psi(\tau)\right)dt\\
&=& \tau \psi'(\tau)-\psi(\tau)\nonumber.
\end{eqnarray}
Hence, $\tau$ should satisfy the equation:
\begin{equation}
\label{taudelta}
\tau \psi'(\tau)-\psi(\tau)=\delta.
\end{equation}
Let $\tau^{*}=\tau(\delta)$ be one solution of \eqref{taudelta}. We use this parameter in order to define the perturbed density modification $f_{x_{i\tau}}(\cdot)$ defined by \eqref{dkl}.
\par 
Now, let us consider the function:
\[G(\tau)=\tau \psi'(\tau)-\psi(\tau)-\delta . \] 
This function has a global minimum at $\tau=0$: $\left.G'(\tau)\right|_{\tau=0}=0, \left.G''(\tau)\right|_{\tau=0}=\psi''(\tau)>0$ and $G(0)=-\delta<0$ for $\delta>0$. Moreover, $G'(\tau)=\tau \psi''(\tau)$: $G'(\tau)<0, (\tau<0)$ and $G'(\tau)>0, (\tau>0)$. Thus, $G(\tau)$  is strictly decreasing for $\tau<0$ and $G(\tau)$ is strictly increasing  for $\tau>0$. Hence, the function $G(\tau)$  has not more than two zeros $\tau_{1}<0$ and $\tau_{2}>0$, if both of them $\tau_{1}$ and $\tau_{2}$ fall into domain of the function $\psi(\tau)$.
\par
For every fixed level of $\delta$, we can study two possible effects of the perturbation \eqref{dkl}. We denote the corresponding perturbed densities by $f_{x_{i\tau_{1}}}$ and $f_{x_{i\tau_{2}}}$. Then, the joint perturbed probability density is expressed as:
\[f_{\bar{x}_{i\tau_{j}}}(\bar{x})=f_{x_{i\tau_{j}}}\prod_{k=1,k\neq i}^{d}f_{x_{k}}(x_{k}), \ \ j=1,2.\]
The corresponding value of the perturbed failure probability $p_{i\delta_{j}} \ \ (j=1,2)$ can be computed as the following integral:
\begin{equation}
\label{jointftau}
p_{i\delta_{j}}=\mathbb{E}_{f_{\bar{x}_{i\tau_{j}}}}\left[\mathbf{I}_{g(\bar{x})<0}\right]=\int\mathbf{I}_{g(\bar{x})<0}f_{\bar{x}_{i\tau_{j}}}d\bar{x},\ \ j=1,2.
\end{equation}
In the same way, the interaction effect can be estimated by perturbing two variables $x_{i}$ and $x_{j}$ at the same time by $\delta_{1}$ and $\delta_{2}$ respectively. Suppose that $\tau(\delta_{1})=(\tau_{1\delta_{1}}, \tau_{2\delta_{1}})$ and $\tau(\delta_{2})=(\tau_{1\delta_{2}}, \tau_{2\delta_{2}})$ are the solutions of equation \eqref{taudelta}, where $\delta_{1}$ and $\delta_{2}$ are the perturbations of KL divergence \eqref{kl} for the variables $x_{i}$ and $x_{j}$, respectively. The new joint probability density function is:
\begin{equation}
f_{\bar{x}_{ij,\tau(\delta_{1}),\tau(\delta_{2})}}(\bar{x})=f_{x_{i\tau(\delta_{1})}}f_{x_{j\tau(\delta_{2})}}\prod_{k=1,k\neq i,j}^{d}f_{x_{k}}(x_{k}).
\end{equation}
The corresponding value of the perturbed failure probability is estimated in the same way by putting in \eqref{jointftau} the new joint probability density $f_{\bar{x}_{ij,\tau(\delta_{1}),\tau(\delta_{2})}}(\cdot)$.
\par
In the next section, we introduce a method to estimate efficiently the perturbed failure probability $p_{i\delta}$ using the same Monte Carlo sample as for the estimation of the original failure probability $p_{f}$. First, we study the effect of the perturbation for different probability distributions.

\subsection{Resulting distributions}
Here, we provide a summary table for the considered distributions from exponential family with the resulting perturbed distribution. For the case when the natural parameter $\eta(\theta)$ is a vector of functions (like normal and log-normal distributions), we propose to analyze the perturbation effect separately for each of the components of  $\eta(\theta)$. For example, for the normal distribution:
\[
\mbox{T}(x)=\left[\begin{array}{c} x\\x^{2}\end{array}\right],\  \ \eta(\mu, \sigma)=\left[\begin{array}{c} \mu/\sigma^{2}\\-1/2\sigma^{2}\end{array}\right].\]
Then we propose to analyze two different density perturbations \eqref{dkl} by taking:
\[
\boldsymbol{\tau}=\left[\begin{array}{c} \tau\\0\end{array}\right] \ \ \text{or}\  \ \boldsymbol{\tau}=\left[\begin{array}{c} 0 \\ \tau\end{array}\right].\]
Table \ref{resdis} provides the results obtained for some distributions from the exponential family: normal, log-normal, exponential and Poisson (as discussed all the calculations are valid for a discrete probability measure as well). As it can be seen, the new perturbed density $f_{x_{i\tau}}$ is still in the same family of distributions  with the perturbed distribution parameters. For the case of normal and log-normal distributions when we study the effect of perturbation of the second component of the natural parameter, there is only one solution for $\tau$: $1-2\tau\sigma^{2}>0$. This solution could not be found analytically but with the help of a numerical solver. For the exponential and the Poisson distributions the solutions for $\tau$ is expressed with the Lambert $W$ function. This stands for the the multivalued inverse relation of the function $f(w)=w\exp(w)$, where $w$ is complex. We denote by $W_{0}(x)$ the upper real branch of the Lambert function on the interval $[-1/\mbox{e},0]$ and by $W_{-1}(x)$ the lower real branch on the same interval.
\par
We will analyze the effect of these perturbations on an analytical example in Section \ref{anfuncex}.
\begin{table}[H]
\centering	 
\resizebox{\textwidth}{!}{
 	\begin{tabular}{c|c|c|c|c|c}
		\hline
Distribution&Natural Parameter&Sufficient Statistics&New cumulant function& $\tau_{1}(\delta)$ and&Resulting distribution\\
&$\eta(\boldsymbol{\theta})$& $\mbox{T}(x)$ & $\psi(\tau)$& $\tau_{2}(\delta)$&\\
\hline
Normal $\mathcal{N}(\mu,\sigma^{2})$& $\left[\begin{array}{c} \mu/\sigma^{2}\\-1/2\sigma^{2}\end{array}\right]$& $\left[\begin{array}{c} x\\x^{2}\end{array}\right]$&
$\left[\begin{array}{c}\mu\tau+\frac{\tau^{2}\sigma^{2}}{2}\\ -\frac{\mu^{2}\tau}{1-2\tau\sigma^2}+\frac{1}{2}\log(1-2\tau\sigma^2) \end{array}\right]$& $\left[\begin{array}{c}\tau_{1,2}=\pm\frac{\sqrt{2\delta}}{\sigma}\\Numerical\ \ solver\end{array}\right]$ & $\left[\begin{array}{c}\mathcal{N}(\mu+\tau\sigma^{2},\sigma^{2})\\ \mathcal{N}(\frac{\mu}{1-2\tau\sigma^2},\frac{\sigma^{2}}{1-2\tau\sigma^2})\end{array}\right]$\\ 
&&&&&\\
LogNormal $\mbox{log}\mathcal{N}(\mu,\sigma^{2})$& $\left[\begin{array}{c} \mu/\sigma^{2}\\-1/2\sigma^{2}\end{array}\right]$& $\left[\begin{array}{c} \log x\\ (\log x)^{2}\end{array}\right]$&
$\left[\begin{array}{c}\mu\tau+\frac{\tau^{2}\sigma^{2}}{2}\\ -\frac{\mu^{2}\tau}{1-2\tau\sigma^2}+\frac{1}{2}\log(1-2\tau\sigma^2) \end{array}\right]$& $\left[\begin{array}{c}\tau_{1,2}=\pm\frac{\sqrt{2\delta}}{\sigma}\\Numerical\ \ solver\end{array}\right]$ & $\left[\begin{array}{c}\log\mathcal{N}(\mu+\tau\sigma^{2},\sigma^{2})\\ \log\mathcal{N}(\frac{\mu}{1-2\tau\sigma^2},\frac{\sigma^{2}}{1-2\tau\sigma^2})\end{array}\right]$\\ 
&&&&&\\
Exponential $\mbox{Exp}(\lambda)$ & $-\lambda$ & $x$ & $\log\left( \frac{\lambda}{\lambda-\tau} \right)$ & $\tau_{1,2}(\delta)=\frac{\lambda\left(W_{-1,0}(-\mbox{e}^{-1-\delta})+1\right)}{W_{-1,0}(-\mbox{e}^{-1-\delta})}$& $\mbox{Exp}(\lambda-\tau)$\\
&&&&&\\
Poisson $\mbox{Pois}(\lambda)$ & $\log(\lambda)$ & $x$ & $\lambda(\mbox{e}^{\tau}-1)$ & $\tau_{1,2}(\delta)=W_{-1,0}\left(-\frac{\lambda-\delta}{\mbox{e}\lambda}\right)+1$ & $\mbox{Pois}(\lambda\exp(\tau))$\\
\hline
\end{tabular}
}
\caption{\label{resdis}Resulting distributions.}
\end{table} 
\newpage
\section{Bounded support distribution}
\begin{figure}[H]
\vspace{-10pt}
\centering
\subfigure[Exponential density modification]{
\includegraphics[width=0.4\textwidth]{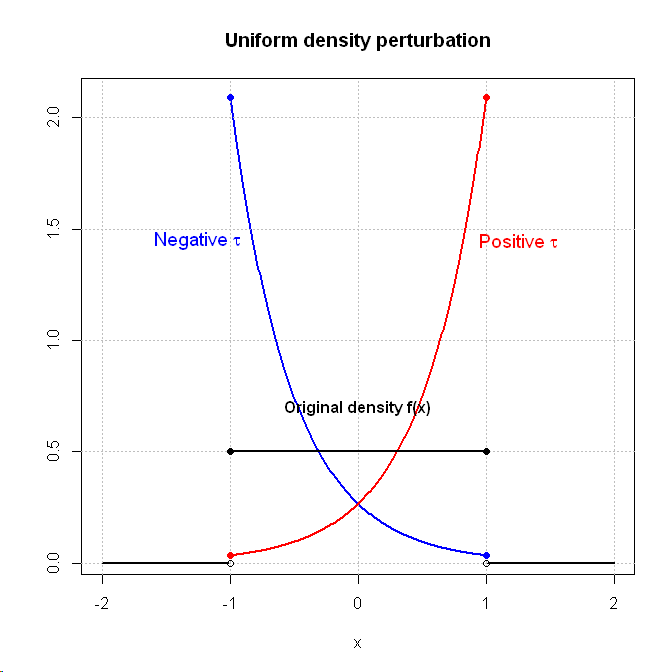}
\label{uniform}}
\subfigure[Boundaries perturbation]{
\includegraphics[width=0.4\textwidth]{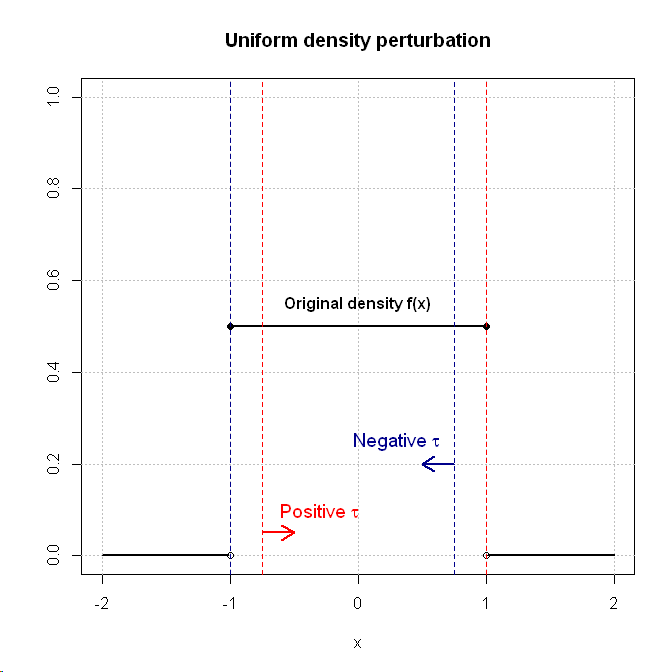}
\label{u_new}}
\caption{Uniform density perturbation.}
\label{uni}
\vspace{-10pt}
\end{figure}

Now, let us consider $x_{i}\thicksim \mbox{U}[a,b]$, the uniform distribution on the interval $[a,b]$. The density is expressed as:
\[ f_{x_{i}}(x)=\frac{1}{b-a}\mathbf{I}_{x\in[a,b]}(x), (b>a).\vspace{-5pt}\]
The uniform distribution does not belong to the family of exponential distributions. It has a limited and a bounded support $[a,b]$. If we apply the same perturbation as for exponential families, the normalization function becomes:
\begin{equation}
\label{tauunif}
\psi(\tau)=\log\left(\frac{\mbox{e}^{\tau b}-\mbox{e}^{\tau a}}{\tau(b-a)}\right), \tau\in\mathbb{R}.\vspace{-5pt}\end{equation}
Then the equation for $\tau(\delta)$ is:\vspace{-5pt}
\[
\frac{\tau b \mbox{e}^{\tau b}-\tau a \mbox{e}^{\tau a}-\mbox{e}^{\tau b}+\mbox{e}^{\tau a}}{\mbox{e}^{\tau b}-\mbox{e}^{\tau a}}-\log\left(\frac{\mbox{e}^{\tau b}-\mbox{e}^{\tau a}}{\tau(b-a)}\right)=\delta.\vspace{-5pt}\]
This equation has no explicit solutions for $\tau$. The solutions can be found using a numeric solver. Suppose that $\tau^{*}=\tau(\delta)$ is a solution of equation \eqref{tauunif}. Then, the perturbed density is:
\[
f_{x_{i\tau}}=\frac{\tau^{*}\mbox{e}^{\tau^{*} x}}{\mbox{e}^{\tau^{*} b}-\mbox{e}^{\tau^{*} a}}\mathbf{I}_{x\in[a,b]}(x).\vspace{-5pt}\]
Therefore, the new perturbed variable $x_{i\tau}$ is no longer uniform on $[a,b]$. This density modification for $a=-1,b=1$ and $\delta=0.5$ is displayed in Figure \ref{uniform}. \par
Notice that working with uncertain variables defined on a compact support, the main source of uncertainty is on the boundaries of the support. For such distributions with a bounded support, we propose to apply another density perturbation. The idea consists in perturbing the original boundaries by $\tau=\pm\delta$. In the same way as with infinite support we consider the effect of positive or negative perturbation. For example, consider $x_{i}\thicksim\mbox{U}[a,b]$ to be uniformly distributed on the interval $[a,b]$. Then, in order to stay inside the support the perturbed random variable $x_{i\tau}$ is uniformly distributed either on $U[a+\delta, b]$ or on $U[a, b-\delta]$. The corresponding density for perturbed uniform distribution can be expressed as:
\[ f_{\tau1}(x)=\frac{1}{b-a-\delta}\mathbf{I}_{x\in[a+\delta,b]}(x) \  \ \text{or}\  \  f_{\tau2}(x)=\frac{1}{b-a-\delta}\mathbf{I}_{x\in[a,b-\delta]}(x) .\]
The same perturbation may be applied to a triangular or a trapezoidal distribution. It can be also applied to the truncated Gaussian distribution if one is interested about the boundary influence on the failure probability $p_{f}$. In this case the density function should be corrected for the new boundaries. 
\noindent
Next, we explain how to estimate efficiently a perturbed failure probability $p_{f_{\delta}}$ with no additional CPU cost.

\section{Inverse importance sampling and sensitivity analysis}
Monte Carlo sampling is one of the most popular simulation methods to estimate a failure probability. We consider the input variables space $\Omega\in\mathbb{R}^{d}$. Recall that all the input variables are independent and that $f_{\bar{x}}(\bar{x})=\prod_{k=1}^{d}f_{x_{k}}(x_{k})$ is the joint density of the input variables. Let $\mathbf{X}^{N}=\{\bar{x}_{1},\ldots,\bar{x}_{N}\}\overset{i.i.d.}{\thicksim} f_{\bar{x}}(\cdot)$ be a sample of size $N$. Then, the estimation of the failure probability $p_{f}$ is given by:
\begin{equation}
\label{pfor}
\widehat{p}_{f}=\frac{1}{N}\sum_{k=1}^{N}\mathbf{I}_{g(\bar{x}_{k})\leq0}.
\end{equation}
Now, assume that $f_{x_{i\tau}}$ is a new perturbed density for the input variable $x_{i}$. Then the new joint density is $f_{\bar{x}_{i\delta}}=f_{x_{1}}\cdots f_{x_{i-1}}f_{x_{i\tau}}f_{x_{i+1}} \cdots f_{x_{d}}= f_{x_{i\tau}}(x_{i})\prod_{k=1,k\neq i}^{d}f_{x_{k}}(x_{k})$. The corresponding failure probability $p_{i\delta}$ is defined as an expectation of the indicator function:
\[
p_{i\delta}=\mathbb{E}_{f_{\bar{x}_{i\delta}}}\left[\mathbf{I}_{g(\bar{x})\leq0}\right]=\int_{\Omega}\mathbf{I}_{g(\bar{x})\leq0}f_{\bar{x}_{i\delta}}(\bar{x})d\bar{x}
\]
Here, we propose to apply the technique used in the Importance Sampling (IS) simulation method. We multiply the integrand function by $\mbox{1}=\frac{f_{\bar{x}}(\bar{x})}{f_{\bar{x}}(\bar{x})}$. Both density functions $f_{\bar{x}}(\bar{x})$ and $f_{\bar{x}_{i\delta}}$ are the products of the density functions of the independent variables $\bar{x}=(x_{1},\ldots,x_{d})\in\Omega\subset\mathbb{R}^{d}$ with the only difference for the variable $x_{i}$. Therefore, we obtain:
\[
p_{i\delta}=\int_{\Omega}\mathbf{I}_{g(\bar{x})\leq0}f_{\bar{x}_{i\delta}}(\bar{x})d\bar{x}=\int_{\Omega}\mathbf{I}_{g(\bar{x})\leq0}\frac{f_{x_{i\tau}}(x_{i})}{f_{x_{i}}(x_{i})}f_{\bar{x}}(\bar{x})d\bar{x}=\mathbb{E}_{f_{\bar{x}}}\left[\mathbf{I}_{g(\bar{x})\leq0}\frac{f_{x_{i\tau}}(x_{i})}{f_{x_{i}}(x_{i})}\right].\]
By doing so, we do not need to throw a new sample according to the unknown density function $f_{\bar{x}_{i\delta}}(\bar{x})$. We are working in the same probability space integrating the function $\left[\mathbf{I}_{g(\bar{x})\leq0}\frac{f_{x_{i\tau}}(x_{i})}{f_{x_{i}}(x_{i})}\right]$. So that, to estimate the perturbed failure probability $p_{i\delta}$ we keep the same sample points from the failure region: $\mathbf{X}_{f}^{N}=\{\bar{x}\in\mathbf{X}^{N}:g(\bar{x})\leq0\}$ that provide non-zero values of the indicator function $\mathbf{I}_{g(\bar{x})\leq0}$. The estimation of the failure probability for the perturbed density is expressed as:
\begin{equation}
\label{pfidelta}
\widehat{p}_{i\delta}=\frac{1}{N}\sum_{k=1}^{N}\mathbf{I}_{g(\bar{x})\leq0}\frac{f_{x_{i\tau}}(\bar{x}_{ki})}{f_{x_{i}}(\bar{x}_{ki})}.
\end{equation}
If we are interested in the interaction effects, we perturb the probability densities for the variables $x_{i}$ and $x_{j}$ simultaneously. Then, the new joint density is:
\[
f_{\bar{x}_{ij,\tau(\delta_{1}),\tau(\delta_{2})}}(\bar{x})=f_{x_{i\tau(\delta_{1})}}(x_{i})f_{x_{j\tau(\delta_{2})}}(x_{j})\prod_{k=1,k\neq i,j}^{d}f_{x_{k}}(x_{k}).
\]
Therefore, in this case the new failure probability can be estimated by:
\[\widehat{p}_{ij,{\delta_{1},\delta_{2}}}=\frac{1}{N}\sum_{k=1}^{N}\mathbf{I}_{g(\bar{x})\leq0}\frac{f_{x_{i\tau(\delta_{1})}}(\bar{x}_{ki})f_{x_{j\tau(\delta_{2})}}(\bar{x}_{kj})}{f_{x_{i}}(\bar{x}_{ki})f_{x_{j}}(\bar{x}_{kj})}.\]
This describes the interaction effects of two variables $x_{i}$ and $x_{j}$ on the failure probability $p_{f}$.
 \par
The proposed reliability sensitivity analysis is based on the analysis of the value of the perturbed failure probability. In order to clearly differentiate the magnitude of the influence to the failure probability, we propose a sensitivity measure for the input variables. We provide the indices formulation in the next section and we study their statistical properties in \ref{indprops}.

\section{Sensitivity indices formulation}
\label{baseindices}
There are many possible choices of sensitivity indices.  In this paper, we propose the sensitivity indices based on the difference $p_{i\delta}-p_{f}$ with the original failure probability $p_{f}$. It is expressed as ratio:
\begin{equation}
\label{sensbas}
S_{i\delta}=\frac{p_{i\delta}-p_{f}}{p_{f}}
\end{equation}

Another possible formulation could be found in Lemaitre and Sergienko et. al, 2012 \citep{paul2012}.\par

The support of $S_{i\delta}$ formulated in \eqref{sensbas} is $[-1,+\infty)$. A negative value of $S_{i\delta}$ means that the proposed density modification reduces the failure probability. Conversely, a positive value of this index means an increase in the failure probability. Zero value of $S_{i\delta}$ means that the variable $x_{i}$ has no impact on the failure probability. \par
In practice, $p_{f}$ is estimated by $\widehat{p}_{f}$ with the Monte Carlo simulation method. In the same way, according to \eqref{pfidelta} $p_{i\delta}$ is estimated by $\widehat{p}_{i\delta}$. The estimator of the indices $\widehat{S}_{i\delta}$ can be expressed as:
\[
\widehat{S}_{i\delta}=\frac{\widehat{p}_{i\delta}-\widehat{p}_{f}}{\widehat{p}_{f}}.
\]
This estimation provides an asymptotic unbiased estimation of $S_{i\delta}$. Moreover, according to the Central Limit Theorem (CLT) and the $\Delta$ - method we have:
\[
\frac{1}{\sqrt{\mbox{VAR}_{f_{\bar{x}}}\left[\widehat{S}_{i\delta}\right]}}\left(\widehat{S}_{i\delta}-S_{i\delta}\right)\overset{N\rightarrow\infty}\rightarrow\mathcal{N}(0,1)\]
The proof and an expression of the asymptotic variance can be found in \ref{appsensind}. 
\par
The above basic indices are computed separately for every fixed value of the density perturbation $\delta$. It allows to study the effect of perturbation by varying the value of $\delta$. It can help in reliability design optimization by adjusting the distribution parameters of the input variables $\bar{x}=(x_{1},\ldots,x_{d})$ in order to achieve the lowest failure probability $p_{f}$.

\section{Analytical function example}
\label{anfuncex}
To investigate the previous indices, we consider a linear function of three independent normally distributed variables: $x_{1},x_{2}, x_{3}\thicksim \mathcal{N}(0,1)$:
\begin{equation}
\label{funcexsens}
g(x_{1}, x_{2}, x_{3})=3-0.1x_{1}-0.5x_{2}-1.0x_{3}.
\end{equation}
As a linear combination of independent Gaussian random variables, $g(\bar{x})$ is distributed as $ \mathcal{N}(\mu, \sigma)$, where $\mu=3$ and $\sigma^{2}=0.1^2+0.5^2+1.0^{2}=1.26$. The original failure probability can be explicitly computed:
\[
p_{f}=1-\Phi(\frac{\mu}{\sigma})=3.7\times10^{-3}.\]
The estimation provided by a Monte Carlo sample of size $N=10^{6}$ yields $\widehat{p}_{f}=3.69\times10^{-3}$. We use the same Monte Carlo sample to estimate the perturbed failure probabilities $p_{i\delta}$ and the sensitivity indices $S_{i\delta}$.\par
\begin{figure}[H]
\vspace{-5pt}
\centering
\subfigure[Negative perturbation]{
\includegraphics[width=0.42\textwidth]{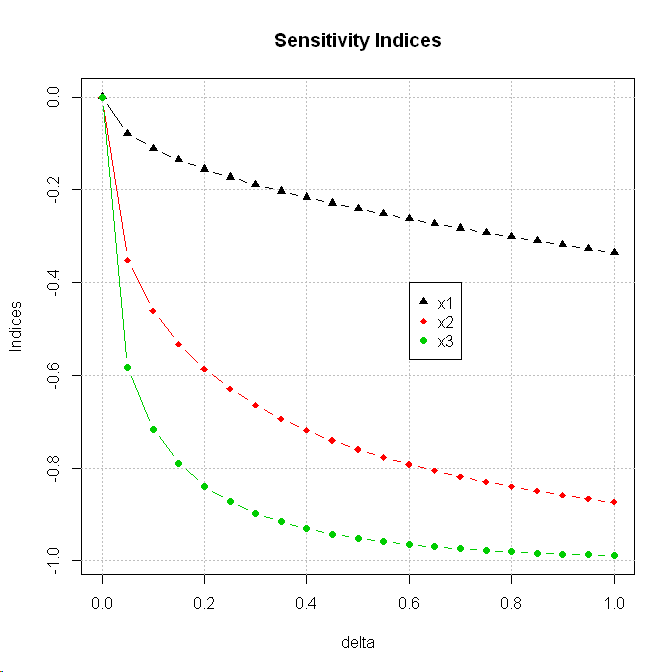}
\label{negtaubas}}
\subfigure[Positive perturbation]{
\includegraphics[width=0.42\textwidth]{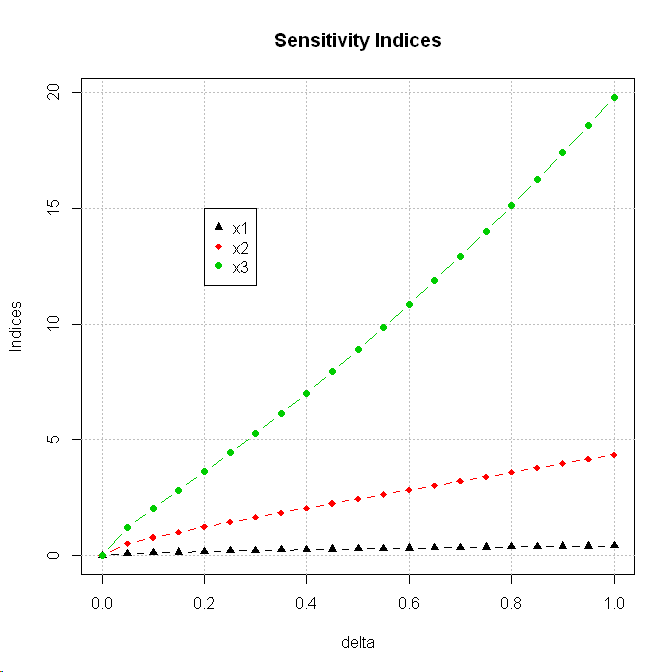}
\label{postaubas}}
\caption{Basic sensitivity indices example. Mean shifting.}
\label{basind}
\vspace{-5pt}
\end{figure}
\begin{figure}[H]
\vspace{-5pt}
\centering
\includegraphics[width=0.42\textwidth]{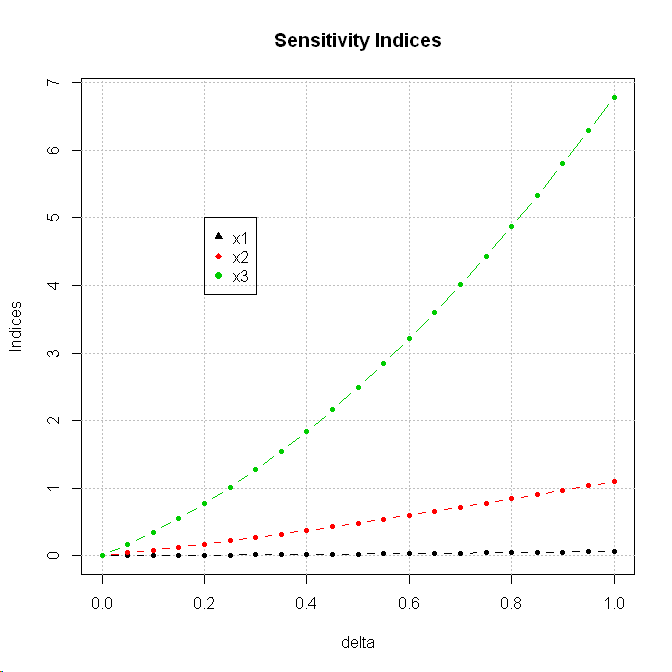}
\caption{Basic sensitivity indices example. Mean and variance shifting.}
\label{anfuncvar}
\vspace{-5pt}
\end{figure}

Recall, that the natural parameter of the normal distribution is given by the vector:\[\eta(\mu, \sigma)=\left[\begin{array}{c}\mu/\sigma^{2}\\-1/2\sigma^{2}\end{array}\right].\]
First, we consider the effect of the perturbation of the first component of natural parameter $\eta$. This perturbation refers as a shifting of the original mean of the distribution (see Table \ref{resdis}). Figure \ref{basind} depicts sensitivity indices calculated for $\delta\in[0,1]$ with the positive \eqref{postaubas} and the negative \eqref{negtaubas} values of perturbation parameter $\tau$. It can be clearly observed that the highest impact on the failure probability is due to the variable $x_{3}$ and that the variable $x_{1}$ has the lowest impact for both cases of $\tau$. Moreover, the higher the value of $\delta$, the higher the influence of the variables.
Now, we study the effect of the perturbation of the second component of the vector of natural parameters $-1/2\sigma^{2}$. In this case, we perturb both the mean and the variance of the original distribution (see Table \ref{resdis}). We use a numerical solver to find the solutions for $\tau$. In this case we consider the only possible effect of perturbation when the normalization function $\psi(\tau)$ is defined and $1-2\tau\sigma^{2}>0$. Figure \ref{anfuncvar} displays the estimated sensitivity indices in this case. We can observe the same relationship: the variable  $x_{3}$ has the highest influence on the failure probability. However, the magnitude of the influence differs from the case where we  only consider a mean shifting. The variable  $x_{1}$ still has almost negligible impact. 

\section{{CO2} storage case example}
In this section, we consider a $\mbox{CO}_{2}$ storage reservoir simulation model. As discussed in the introduction, subsurface $\mbox{CO}_{2}$ storage is always associated with an excess reservoir pressure. One of the primary environmental risks is a pressure-driven leakage of  $\mbox{CO}_{2}$ from the storage formation.\par
\begin{figure}[H]
\vspace{-5pt}
\centering
\includegraphics[width=0.4\textwidth]{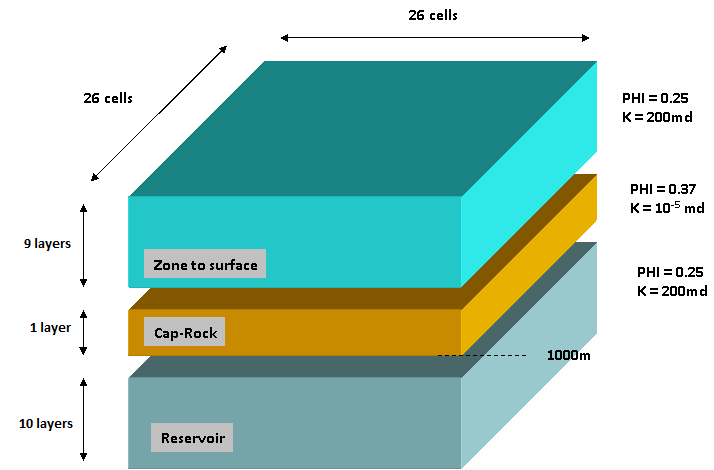}
\caption{\label{co2model}Reservoir model.}
\vspace{-10pt}
\end{figure}

In order to assess the risk of  $\mbox{CO}_{2}$ leakage through the cap rock we consider a synthetic reservoir model. The structure of the reservoir is reduced to its simplest expression. The model is made up of three zones (Figure \eqref{co2model}): a reservoir made of 10 layers, a cap-rock made up of 1 layer,  a zone-to-surface composed of 9 layers.
\noindent
The XY size of the grid is set at 10 km total length representing 26x26x20 model grid. Each layer is 5m thick, including the cell above the cap-rock. 	The zone above the cap-rock (up to the surface) is currently set to 1 layer. The salinity of the water is 35gm/l. The temperature of the reservoir is set to 60C and the initial pressure is hydrostatic. The injection bottom rate is set to $10^{6}$ tons/year. The fracture pressure is estimated by geomechanical experts to be $\mbox{P}_{fracture}=122$ bars. Exceeding this value during the injection can lead to a leakage. The simulation period is 55 years that include an injection period of 15 years followed by 40 years of storage. In this study we analyze the possibility of a leakage through a cap rock. Therefore, we consider pressure in the storage reservoir at the last year of injection as an objective function.
\par
The uncertain variables selected for this study characterize the reservoir and the fluid properties. It implies different $\mbox{CO}_{2}$ flowing possibilities between the reservoir layers. Table \eqref{tab:UP} represents the variables description with their range of minimum and maximum values. 
\begin{table}[H]
	 \centering
 	\begin{tabular}{c|l|cc}
		\hline
				Name&Description&Min&Max\\
		\hline
	PORO&Reservoir Porosity&0.15&0.35\\
	KSAND&Reservoir Permeability&10&300\\
	KRSAND&Water relative permeability end-point&0.5&1.0\\
\hline
		\end{tabular}
		\caption{\label{tab:UP}Uncertain variables.}
\end{table} \noindent
For sake of clarity we transform the original intervals into $[-1,1]$. In this section, we assume the truncated standard normal distribution for all the input variables $\mathcal{N}_{[-1,1]}(0,1)$.
\par
The performance function for this example can be formulated. Suppose that $\mbox{P}_{reservoir}(\bar{x})$ is a function of the reservoir pressure depending on the input variables configuration $\bar{x}\in\Omega\subset\mathbb{R}^{d}$. Then, the performance function defining the event of a gas leakage can be expressed as:
\[ g(\bar{x})=\mbox{P}_{fracture}-\mbox{P}_{reservoir}(\bar{x}).\]
The reservoir pressure $\mbox{P}_{reservoir}(\bar{x})$ is computed with a complex dynamic reservoir simulator. For this reason, we use the Gaussian process based response surface model approximation $\widehat{P}(\bar{x})$ \citep{krig,Sacks1989,Welch1992,Santner2003}. By approximating the function of the reservoir pressure, we can quantify the risk and estimate the reliability of the system. 
\par The original reference value of the failure probability computed by GP model approximation with  Monte Carlo sample of size $N=10^{6}$ provides an estimation: $\widehat{p}_{f}=2.26\times 10^{-4}$. We keep this sample to estimate the perturbed failure probability $p_{i\delta}$ and the sensitivity indices $S_{i\delta}$. 
\par
\begin{figure}[H]
\vspace{-10pt}
\centering
\subfigure[Negative perturbation]{
\includegraphics[width=0.42\textwidth]{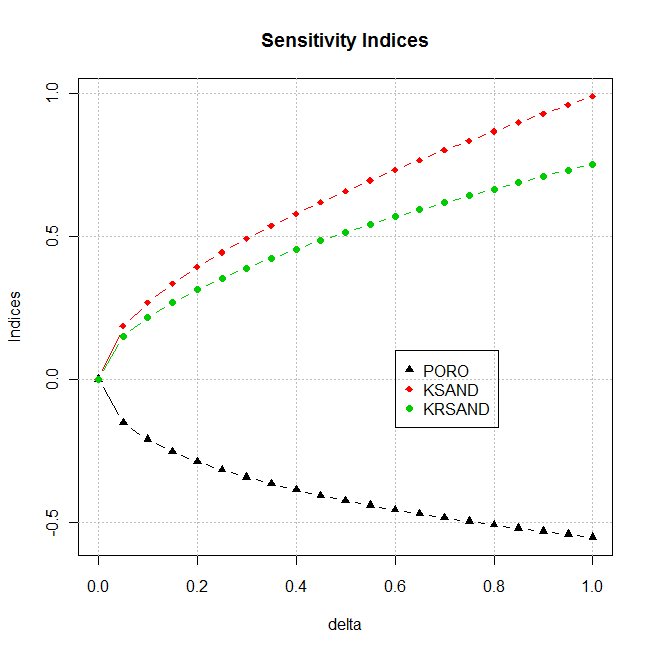}
\label{co2exneg}}
\subfigure[Positive perturbation]{
\includegraphics[width=0.42\textwidth]{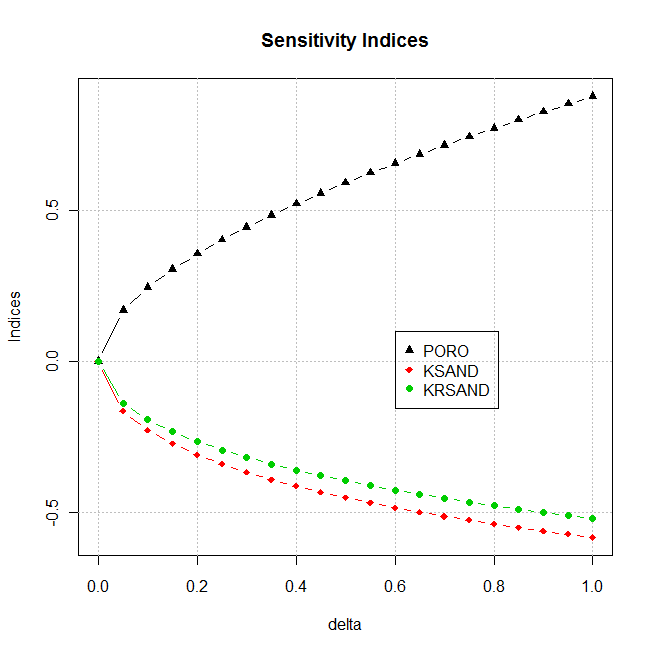}
\label{co2expos}}
\caption{Mean shifting.}
\label{co2exmean}
\vspace{-10pt}
\end{figure}

There are two possible ways to perturb the truncated Gaussian distribution. We can use the perturbation defined for the exponential family or we can study the effect of the perturbation of the distribution support boundaries. In this example, we compare the results for both cases. We start with the sensitivity indices calculated by the mean shifting. Figure \ref{co2exmean} displays the evolution of the sensitivity indices for $\delta\in[0,1]$ for negative (Figure \ref{co2exneg}) and positive (Figure \ref{co2expos}) values of shifting. The variables ranking depends on the sign of $\tau$. When $\tau>0$ (i.e. the positive mean shifting) the porosity variable has the highest impact on the failure probability. It means that increasing the mean value of the reservoir porosity $\mbox{PORO}$ leads to increasing the failure probability. On the contrary, increasing the mean value of the reservoir permeability $\mbox{KSAND}$ and the water relative permeability end-point $\mbox{KRSAND}$ has a negative effect on the risk of leakage. For the negative mean shifting, the variables $\mbox{KSAND}$ and $\mbox{KRSAND}$ have the highest influence on the failure probability. Reducing the reservoir permeability and the end-point water relative permeability impedes the gas flow in the reservoir. It increases the risk of the leakage.
\par
\begin{figure}[H]
\vspace{-10pt}
\centering
\subfigure[Negative shifting]{
\includegraphics[width=0.42\textwidth]{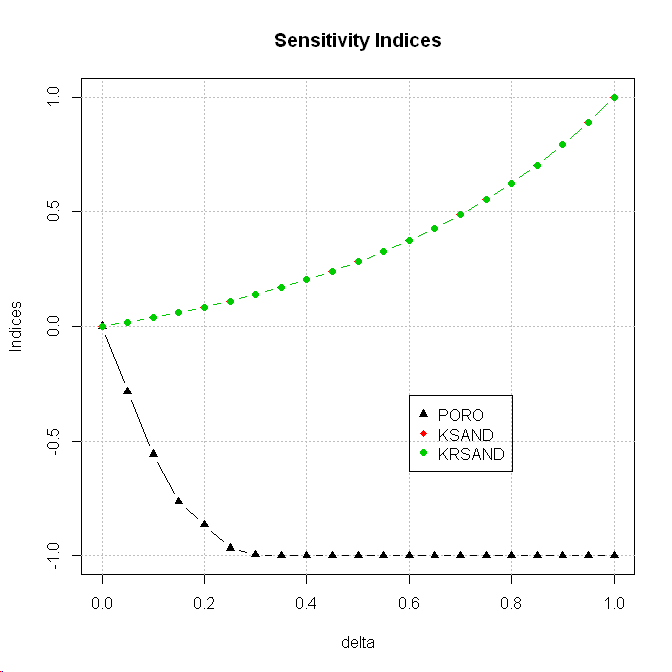}
\label{co2exnegbo}}
\subfigure[Positive shifting]{
\includegraphics[width=0.42\textwidth]{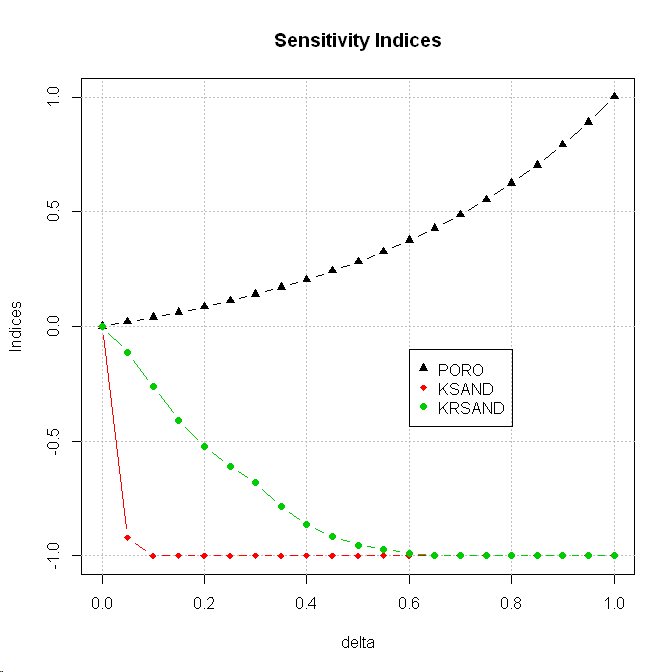}
\label{co2exposbo}}
\caption{Boundaries shifting.}
\label{co2exnewbo}
\vspace{-10pt}
\end{figure}

Now, we consider the boundaries perturbation. It means that we are moving one of the the distribution boundaries by $\delta$ in positive or negative directions keeping the values of mean and variance unchanged. Figure \ref{co2exnewbo} depicts the sensitivity indices for $\delta\in[0,1]$ for negative (Figure \ref{co2exnegbo}) and positive (Figure \ref{co2exposbo}) values of $\tau=\pm\delta$. 
\noindent When $\tau<0$ (i.e. the resulting distribution is $\mathcal{N}_{[-1,1-\delta]}(0,1)$) the porosity variable has the highest impact on the failure probability. It means that by decreasing the maximum value of the reservoir porosity $\mbox{PORO}$ the failure probability decreases.  It is also shown by Figure \ref{co2exneg}. It can be also observed that for this variable the new perturbed failure probability $p_{i\delta}$ is equal to zero when $\delta>0.3$.  When $\tau>0$ the resulting distribution is $\mathcal{N}_{[-1+\delta,1]}(0,1)$. For this case, increasing the reservoir permeability $\mbox{KSAND}$ and the water relative permeability end-point $\mbox{KRSAND}$ reduces the failure probability $p_{i\delta}$. 
\par
Both methods provide comprehensive and complementary results. If the main uncertainty is about the boundaries the one can start with the boundaries perturbation. By moving the boundaries of the original distribution, the one can determine the safe intervals for the input variables by detecting the value of $\delta$: $p_{i\delta}=0$ or $S_{i\delta}=-1$. After that, the effect of the mean perturbation can be studied.

\section{Conclusions}
In this paper, we have studied and adapted a recently introduced approach to the reliability sensitivity analysis. Currently the majority of the methods for reliability analysis is based on the variance decomposition and Sobol' sensitivity indices. We present a moment independent sensitivity measure. The method is based on a perturbation of the original probability distribution of the input random variables. We can analyze the \textit{a priori }assumption about the input distributions and measure the effect of some possible deviations from this assumption. In particular, we select the Kullback-Leibler divergence as a measure of the perturbation.
\par
We have provided different possible density perturbations with the resulting distributions for the exponential family of distributions. We have also studied the distributions with a bounded support and the effect of the boundaries perturbation.
\par
Considering a proposed perturbation for an input variable, we present an effective method to estimate the corresponding perturbed failure probability. The method is based on a technique coming from the importance sampling simulation method. It allows to estimate the new failure probability without supplementary performance function evaluations. The new sensitivity indices formulation describes the relationship between the new failure probability $p_{i\delta}$ and the original failure probability of the system $p_{f}$. By varying the value of the perturbation $\delta$, we can study how the positive or negative probability density perturbation affects the failure probability. If the model has controllable input variables, the method can help improving the system reliability and the design optimization. The presented analysis on the statistical properties of the proposed estimators for the perturbed failure probability $\widehat{p}_{i\delta}$ and the sensitivity indices $\widehat{S}_{i\delta}$ shows asymptotic normality of the estimators.
\par
We investigated the method on an analytical and a $\mbox{CO}_{2}$ storage reservoir cases. The method provides promising results and can be applied in the reliability sensitivity analysis. 
\section{Acknowledgments}
This work has been partially supported by the French National Research Agency (ANR) through COSINUS program (project COSTA-BRAVA n ANR-09-COSI-015).
The authors also would like to thank Nicolas Maurand and Dan Bossie-Codreanu for the presented $\mbox{CO}_{2}$ reservoir model. We
thank S\'ebastien Da-Veiga, Nicolas Bousquet and Bertrand Iooss for fruitful discussions.

\appendix
\section{Properties of the normalization function }
\label{psiapp}
For a given random  variable $x$ with probability density $f(x)$. We propose the density modification $x_{\tau}\thicksim f_{\tau}(\cdot)$ as follows:
\vspace{-5pt}
\begin{equation*}
f_{\tau}(x)=\exp(\tau x-\psi(\tau))f(x),
\vspace{-5pt}
\end{equation*}
where $\psi(\tau)$ is a normalization function given by:
\vspace{-5pt}
\[ \psi(\tau)=\log\left[\int_{-\infty}^{\infty}\exp(\tau x)f(x)dx\right].\vspace{-5pt}\]
Let us define  $\mathcal{D}=\{\tau\in\mathbb{R}:\psi(\tau)<+\infty\}$. $\mathring{\mathcal{D}}$ defines the interior of $\mathcal{D}$.  We will also suppose, that $\exists\varepsilon: \mathcal{D}\supset]-\varepsilon,\varepsilon[$. Here, we will study the properties of this normalization function.

\subsection{Derivatives}
\begin{itemize}
	\item $\psi'(\tau)=\mathbb{E}[x_{\tau}], \tau\in \mathring{\mathcal{D}}$
\end{itemize}\vspace{-5pt}
\[ \psi'(\tau)=\frac{\frac{d}{d\tau}\left[\int_{-\infty}^{\infty}\exp(\tau x)f(x)dx\right]}{\int_{-\infty}^{\infty}\exp(\tau x)f(x)dx}=\int_{-\infty}^{\infty}x\exp(\tau x-\psi(\tau))f(x)dx=\int_{-\infty}^{\infty}x f_{\tau}(x)dx=\mathbb{E}[x_{\tau}]\vspace{-5pt}\]
\begin{itemize}
	\item $\psi''(\tau)=\mathbb{E}\left[x_{\tau}-\mathbb{E}(x_{\tau})\right]^{2}, \tau\in \mathring{\mathcal{D}}$
\end{itemize}
\vspace{-5pt}
\begin{eqnarray*}
\psi''(\tau)&=&\frac{d}{d\tau}\left[\int_{-\infty}^{\infty} x\exp(\tau x -\psi(\tau)) f(x)dx \right]=\\
&=&\int_{-\infty}^{\infty}\left[x(x-\psi'(\tau))\exp(\tau x -\psi(\tau)) f(x)dx\right]=\\
&=&\int_{-\infty}^{\infty}x^{2}f_{\tau}(x)dx-\psi'(\tau)\int_{-\infty}^{\infty}x f_{\tau}(x)dx=\int_{-\infty}^{\infty}x^{2}f_{\tau}(x)dx-\left[\mathbb{E}[x_{\tau}]\right]^{2}=\\
&=&\mathbb{E}\left[x_{\tau}-\mathbb{E}(x_{\tau})\right]^{2}=\mbox{VAR}(x_{\tau})
\vspace{-5pt}
\end{eqnarray*}
\begin{itemize}
	\item $\psi'''(\tau)=\mathbb{E}\left[x_{\tau}-\mathbb{E}(x_{\tau})\right]^{3}, \tau\in \mathring{\mathcal{D}}$
\end{itemize} 
\vspace{-15pt}
\begin{eqnarray*} 
\mathbb{E}\left[x_{\tau} - \mathbb{E}(x_{\tau})\right]^{3}&=&\int_{-\infty}^{\infty}\left[x-\psi'(\tau)\right]^{3}f_{\tau}(x)dx=\\
&=&\int_{-\infty}^{\infty}x^{3}f_{\tau}(x)dx-3\psi'(\tau)\int_{-\infty}^{\infty}x^{2}f_{\tau}(x)dx+3\psi'(\tau)^2\int_{-\infty}^{\infty}x f_{\tau}(x)dx-\psi'(\tau)^{3}\\
&=&\int_{-\infty}^{\infty}x^{3}f_{\tau}(x)dx-3\psi'(\tau)\int_{-\infty}^{\infty}x^{2}f_{\tau}(x)dx+2\psi'(\tau)^{3}=\psi'''(\tau)
\vspace{-5pt}\end{eqnarray*}

\section{Statistical properties of the indices estimator}
 \label{indprops}
Here, we will study statistical properties of the estimator of the perturbed failure probability $\widehat{p}_{i\delta}$ and corresponding estimator of the sensitivity indices $S_{i\delta}=\frac{\widehat{p}_{i\delta}-\widehat{p}_{f}}{\widehat{p}_{f}}$. We will start with studying the properties of  $\widehat{p}_{i\delta}$.
\subsection{Estimator of the perturbed failure probability}
Suppose, $f_{\mathbf{x}}(\mathbf{x})=\prod_{i=1}^{d}f_{x_{i}}(x_{i})$ is the input joint density and $f_{x_{i\tau}}$ is a perturbed probability density for the variable $x_{i}$. Recall that for a sample of size $N$: $\{\mathbf{x}_{1},\ldots,\mathbf{x}_{N}\}\overset{i.i.d.}{\thicksim}f_{\mathbf{x}}(\mathbf{x})$, the estimation of $\widehat{p}_{i\delta}$ is computed by: \vspace{-5pt}
\[ \widehat{p}_{i\delta}=\frac{1}{N}\sum_{k=1}^{N}\mathbf{I}_{g(\mathbf{x}_{k}<0)}\frac{f_{x_{i\tau}}(\mathbf{x}_{k_{i}})}{f_{x_{i}}(\mathbf{x}_{k_{i}})}. \vspace{-5pt}\]
First, we study the expectation and the variance of this estimator.
\begin{enumerate}
	\item $\mathbb{E}_{f_{\mathbf{x}}}\left[ \widehat{p}_{i\delta} \right]=p_{i\delta}$
	\item $\mbox{VAR}_{f_{\mathbf{x}}}\left[ \widehat{p}_{i\delta} \right]=\frac{1}{N}\mbox{VAR}_{f_{\mathbf{x}}}\left[\mathbf{I}_{g(\mathbf{x})<0}\frac{f_{x_{i\tau}}(x_{i})}{f_{x_{i}}(x_{i})}\right]=\frac{1}{N}\left[\int \mathbf{I}_{g(\mathbf{x})<0}\frac{f^{2}_{x_{i\tau}}(x_{i})}{f^{2}_{x_{i}}(x_{i})}f_{\mathbf{x}}(\mathbf{x})d\mathbf{x}-p^{2}_{i\delta}\right]$
\end{enumerate}
This variance tends to $0$ when $N\rightarrow\infty$. Furthermore, by the Central Limit Theorem (CLT):\vspace{-5pt}
\[ \frac{1}{\sqrt{\mbox{VAR}_{f_{\mathbf{x}}}\left[ \widehat{p}_{i\delta} \right]}}\left(\widehat{p}_{i\delta}-p_{i\delta}\right)\overset{N\rightarrow\infty}\rightarrow\mathcal{N}(0,1) \vspace{-5pt}\]

Note that the covariance between the estimator $\widehat{p}_{f}$ and $\widehat{p}_{i\delta}$ does not vanish. Indeed, we use the same sample to estimate $p_{f}$ and ${p}_{i\delta}$. We can compute this covariance:
\begin{eqnarray*}
\mbox{COV}(\widehat{p_{f}},\widehat{p}_{i\delta})&=&\mathbb{E}_{f_{\mathbf{x}}}(\widehat{p_{f}}\widehat{p}_{i\delta})-\mathbb{E}(\widehat{p_{f}})\mathbb{E}(\widehat{p}_{i\delta})=\\
&=&\frac{N}{N^{2}}\int\mathbf{I}_{g(\mathbf{x}_{k}<0)}\frac{f_{x_{i\tau}}(x_{i})}{f_{x_{i}}(x_{i})}{f_{x_{i}}(x_{i})}d\mathbf{x}-p_{f}{p}_{i\delta}=\\
&=&\frac{1}{N}p_{i\delta}(1-p_{f})\end{eqnarray*}
The value of this covariance decreases when the sample size  $N$ increases.  

\subsection{Sensitivity indices estimator}
\label{appsensind}
Recall, that the first sensitivity index is:
\[ S_{i\delta}=\frac{p_{i\delta}-p_{f}}{p_{f}}=\frac{p_{i\delta}}{p_{f}}-1, i=1,\ldots,d.\vspace{-5pt} \]
We estimate this value with the Monte Carlo method and the importance sampling by estimating consistently $p_{i\delta}$  and $p_{f}.$   The estimator of this index  is:\vspace{-5pt}
\[ \widehat{S}_{i\delta}=\frac{\widehat{p}_{i\delta}}{\widehat{p}_{f}}-1.\vspace{-5pt}\]
Here, we will study some proprieties of this estimator. 
\par
It is not straightforward to compute directly $\mathbb{E}_{f_{\mathbf{x}}}\left[\widehat{S}_{i\delta}\right]=\mathbb{E}_{f_{\mathbf{x}}}\left[\frac{\widehat{p}_{i\delta}}{\widehat{p}_{f}}\right]-1$ and $\mbox{VAR}_{f_{\mathbf{x}}}\left[\widehat{S}_{i\delta}\right]=\mbox{VAR}_{f_{\mathbf{x}}}\left[\frac{\widehat{p}_{i\delta}}{\widehat{p}_{f}}\right]$.
We propose to use the Delta Method to approximate these values \citep{van2000asymptotic}.
\par
Let us recall the Taylor expansion with integral form for the remainder. Let $\phi$ be a two times differentiable function on $[t_{0},t]$, then:
\begin{equation}
\label{taylor}
\phi(t):=\phi(t_{0})+\phi'(t_{0})(t-t_{0})+\int_{t_{0}}^{t}(1-u)\phi''(u)du
\end{equation}
We will define a function: \vspace{-5pt}
\[\phi(t)=\frac{y(t)}{x(t)}-1,\vspace{-5pt}\]
where $x(t)=(1-t)p_{f}+t\widehat{p}_{f}$ and $y(t)=(1-t)p_{i\delta}+t\widehat{p}_{i\delta}$. For this function: $\phi(0)=S_{i\delta}$ and $\phi(1)=\widehat{S}_{i\delta}$.  Following the Taylor expansion \eqref{taylor} we will expand $\phi(t)$ with $t=1$ and $t_{0}=0$. First, we will compute the derivatives. 
\begin{enumerate}
	\item $x'(t)=\widehat{p}_{f}-p_{f}$
	\item  $y'(t)=\widehat{p}_{i\delta}-p_{i\delta}$
\end{enumerate}
Then,\vspace{-15pt}
\begin{eqnarray*}
\phi'(t)&=&\frac{y'(t)x(t)-x'(t)y(t)}{x^{2}(t)}=\frac{p_{f}\widehat{p}_{i\delta}-\widehat{p}_{f}p_{i\delta}}{\left((1-t)p_{f}+t\widehat{p}_{f}\right)^{2}}\\
\left.\phi'(t)\right|_{t=0}&=&\frac{p_{f}\widehat{p}_{i\delta}-\widehat{p}_{f}p_{i\delta}}{p^{2}_{f}}
\end{eqnarray*} 
The second derivative is:\vspace{-5pt}
\[ \phi''(t)=\frac{2x(t)x'(t)\left({\widehat{p}_{f}p_{i\delta}}-p_{f}\widehat{p}_{i\delta}\right)}{x^{4}(t)}=\frac{2x'(t)\left({\widehat{p}_{f}p_{i\delta}}-p_{f}\widehat{p}_{i\delta}\right)}{x^{3}(t)}=\frac{2(\widehat{p}_{f}-p_{f})\left({\widehat{p}_{f}p_{i\delta}}-p_{f}\widehat{p}_{i\delta}\right)}{x^{3}(t)}.\vspace{-5pt}\]
Therefore, the reminder is:\vspace{-5pt}
\[ \int_{0}^{1}\frac{2(1-t)(\widehat{p}_{f}-p_{f})\left({\widehat{p}_{f}p_{i\delta}}-p_{f}\widehat{p}_{i\delta}\right)}{\left((1-t)p_{f}+t\widehat{p}_{f}\right)^{3}}dt=\frac{(\widehat{p}_{f}-p_{f})(p_{f}\widehat{p}_{i\delta}-\widehat{p}_{f}p_{i\delta})}{p^{2}_{f}\widehat{p}_{f}}.\vspace{-5pt}\]
So that, by the Taylor expansion \eqref{taylor} we obtain:\vspace{-5pt}
\[ \widehat{S}_{i\delta}=S_{i\delta}+\frac{p_{f}\widehat{p}_{i\delta}-\widehat{p}_{f}p_{i\delta}}{p^{2}_{f}}+\frac{(\widehat{p}_{f}-p_{f})(p_{f}\widehat{p}_{i\delta}-\widehat{p}_{f}p_{i\delta})}{p^{2}_{f}\widehat{p}_{f}}.\vspace{-5pt}\]
The last term $R=\frac{(\widehat{p}_{f}-p_{f})(p_{f}\widehat{p}_{i\delta}-\widehat{p}_{f}p_{i\delta})}{p^{2}_{f}\widehat{p}_{f}}$ is the remainder. This remainder is bounded and we can neglect it for the approximation.\vspace{-5pt}
\[ \widehat{S}_{i\delta}\approx S_{i\delta}+\frac{p_{f}\widehat{p}_{i\delta}-\widehat{p}_{f}p_{i\delta}}{p^{2}_{f}}.\vspace{-5pt} \]
Now, we can approximate the mean and the variance of $\widehat{S}_{i\delta}$.
\begin{spacing}{1.25}
\begin{enumerate}
\item $\mathbb{E}_{f_{\mathbf{x}}}\left[\widehat{S}_{i\delta}\right]\thicksim S_{i\delta}$
\item $\mbox{VAR}_{f_{\mathbf{x}}}\left[\widehat{S}_{i\delta}\right]\thicksim  \mbox{VAR}_{f_{\mathbf{x}}}\left[ \frac{\widehat{p}_{i\delta}}{p
_{f}}\right]+ \mbox{VAR}_{f_{\mathbf{x}}}\left[ \frac{\widehat{p}_{f}p_{i\delta}}{p^{2}_{f}}\right]-\frac{2p_{i\delta}}{p_{f}^{3}}\mbox{COV}(\widehat{p}_{i\delta},\widehat{p}_{f})=
\frac{1}{p_{f}^{2}} \mbox{VAR}_{f_{\mathbf{x}}}\left[ \widehat{p}_{i\delta} \right]-\frac{p^{2}_{i\delta}(1-p_{f}) } {N p_{f}^{3}}.$
\end{enumerate}\end{spacing} \noindent
Therefore, the variance of the indices estimator tends to $0$ when $N\rightarrow\infty$. With some extra computations we can show that:\vspace{-5pt}
\[ \frac{1}{\sqrt{\mbox{VAR}_{f_{\mathbf{x}}}\left[\widehat{S}_{i\delta}\right]}}\left(\widehat{S}_{i\delta}-S_{i\delta}\right)\overset{N\rightarrow\infty}\rightarrow\mathcal{N}(0,1).\vspace{-5pt}\]
Knowing the variance of the estimator, the confidence region for  the indices  may be computed.

\bibliographystyle{elsarticle-num}
\bibliography{references}
\end{document}